\documentclass[english]{article}
\usepackage[T1]{fontenc}
\usepackage[latin9]{inputenc}
\usepackage{amsthm}
\usepackage{amsmath}
\usepackage{amssymb}
\usepackage{esint}
\DeclareMathOperator{\Id}{Id}

\makeatletter
\theoremstyle{plain}
\newtheorem{thm}{Theorem}
  \theoremstyle{definition}
  \newtheorem{defn}[thm]{Definition}

\usepackage{mathrsfs}

\makeatother
\usepackage{babel}

\begin{document}

\title{Some generic properties of non degeneracy for critical points of
functionals and applications}

\author{M. Ghimenti,\thanks{Dipartimento di Matematica Applicata ``U.Dini'', via Buonarroti 1c, 56127, Pisa, Italy.
email:{\tt marco.ghimenti@dma.unipi.it, a.micheletti@dma.unipi.it }}
 A.M. Micheletti\addtocounter{footnote}{-1}
\footnotemark}
\date{}	
\maketitle

\section{Introduction}

In these last years there have been some theorems of existence and
multiplicity of solutions of the following equation\begin{equation}
\left\{ \begin{array}{ccc}
-\varepsilon^{2}\Delta_{g}u+u=|u|^{p-2}u &  & \text{in }M\\
u\in H_{g}^{1}(M)\end{array}\right..\label{eq:main}\end{equation}
Here $(M,g)$ is a smooth connected compact Riemannian manifold of
dimension $n\ge3$ embedded in $\mathbb{R}^{N}$. In \cite{BBM07,H09,V08}
it is shown that the number of solutions is influenced by the topology
of $M$. In \cite{DMP09,MP09b,MP09} there are some results about
the effect of the geometry of $M$ in finding solutions, more precisely
the role of the scalar curvature $S_{g}$ of $(M,g)$. In these results,
a type of nondegeneracy on the critical points of $S_{g}$ is assumed.
In Section \ref{sec:2} we give a generic property (see Theorem \ref{thm:2-1}
and \cite{MP10}) of nondegeneracy of critical points of $S_{g}$ with
respect to the metric $g$ and an application to the results of the
papers \cite{DMP09,MP09b,MP09} to obtain some theorems of existence
and multiplicity of solutions of (\ref{eq:main}).

In Section \ref{sec:3} we consider the Neumann problem\begin{equation}
\left\{ \begin{array}{ccc}
-\varepsilon^{2}\Delta_{g}u+u=|u|^{p-2}u &  & \text{in }\Omega\\
\frac{\partial u}{\partial\nu}=0 &  & \text{on }\partial\Omega\end{array}\right..\label{eq:neu}\end{equation}

This problem has similar feature with problem (\ref{eq:main}). Indeed
there are some theorems about the existence of solutions of (\ref{eq:neu})
in which the mean curvature of the boundary $\partial\Omega$ of the
domain $\Omega$ plays the same role of the scalar curvature $S_{g}$
of the manifold $(M,g)$ with respect to problem (\ref{eq:main}).
We show (see Theorem \ref{thm:3-1} and \cite{MPpre}) a generic
property for critical points of the mean curvature of the boundary
$\partial\Omega$ with respect to the deformation of the domain $\Omega$.
Thus the results in \cite{DY,DFW,G,GWW,Li,W1,WW} can be applied. 

In Section \ref{sec:4} we consider a Riemannian manifold $(M,g)$
embedded in $\mathbb{R}^{N}$ invariant with respect to a given involution
of $\mathbb{R}^{N}$ and we show some results (see theorems \ref{thm:4-1}
and \ref{thm:4-2} and \cite{GMp}) of genericity for non degenerate sign changing
solutions of the problem (\ref{eq:ptau}) in which the metric $g$
is considered as a parameter. By these results we can use the Morse
theory and we can give an estimate of the number of solutions which
change sign exactly once. To obtain our generic properties of non
degeneracy of critical point the main tool is an abstract transversality
theorem of of Quinn \cite{Qu70}, Uhlenbeck \cite{Uh76}, Saut and
Temam \cite{ST79} which will be recalled in the appendix.

\section{Genericity of nondegeneracy for critical points of the scalar curvature
for a Riemannian manifold and applications\label{sec:2}}

Let $M$ be a connected compact $C^{\infty}$ manifold of dimension
$N\geq2$, without boundary. Let $\mathscr{M}^{k}$ be the set of
all $C^{k}$ Riemannian metrics on $M$. Any $g\in\mathscr{M}^{k}$
determines the scalar curvature $S_{g}$ of $(M.g)$. Our goal is
to prove that for a generic Riemannian metric $g$ the critical points
of the scalar curvature $S_{g}$ are nondegenerate. More precisely
we can prove the following result (see \cite{MP10})
\begin{thm}
The set\[
A=\left\{ g\in\mathscr{M}^{k}\ :\ \text{all the critical points of }S_{g}\text{ are nondegenerate}\right\} \]
is an open dense subset of $\mathscr{M}^{k}$. Here $k\ge3$\label{thm:2-1}.
\end{thm}
In the following we denote by $\mathscr{S}^{k}$ the space of all
$C^{k}$ symmetric covariant $2$-tensors on $M$. $\mathscr{S}^{k}$
is a Banach space equipped with the norm $\|\cdot\|_{k}$ defined
in the following way. We fix a finite covering $\left\{ V_{\alpha}\right\} _{\alpha\in L}$
of $M$ such that the closure of any $V_{\alpha}$ is contained by
$U_{\alpha}$, where $\left(U_{\alpha},\psi_{\alpha}\right)$ is an
open coordinate neighborhood. If $h\in\mathscr{S}^{k}$, denoting
$h_{ij}$ the components of $h$ with respect to local coordinates
$(x_{1},\dots,x_{n})$ on $V_{\alpha}$, we define
\begin{equation}
\|h\|_{k}=
\sum_{\alpha\in L}\ \sum_{|\beta|\leq k}\  \sum_{i,j=1}^n\ 
\sup_{\psi_{\alpha}(V_{\alpha})}\left|\frac{\partial^{\beta}h_{ij}}
{\partial x_{1}^{\beta_{1}}\cdots\partial x_{n}^{\beta_{n}}}\right|.
\label{eq:normak}
\end{equation}
The set $\mathscr{M}^{k}$ of all $C^{k}$ Riemannian metrics on $M$
is an open subset of $\mathscr{S}^{k}$. 

We can apply Theorem \ref{thm:2-1} to study the following problem
\begin{equation}
\left\{ \begin{array}{ccc}
-\varepsilon^{2}\Delta_{g}u+u=|u|^{p-2}u &  & \text{in }M\\
u\in H_{g}^{1}(M)\end{array}\right.\label{eq:main-sec2}\end{equation}
where $p>2$ if $N=2$, $2<p<\frac{2N}{N-2}$ if $N\ge3$. Here $H_{g}^{1}(M)$
is the completion of $C^{\infty}$ with respect to \[
\|u\|_{g}^{2}=\int_{M}(|\nabla_{g}u|^{2}+u^{2})d\mu_{g}.\]

There are some recent results on the effect of the geometry of $(M,g)$
on the number of solutions of (\ref{eq:main-sec2}). We will see that
the scalar curvature $S_{g}$ relative to the metric $g$ is the geometric
property which influences the number of solutions. Indeed we have
the following results about the role of the scalar curvature.
\begin{thm}
\label{thm:2-2}{\rm(See\cite{MP09})} For any $C^{1}$-stable critical
point $\bar{q}$ of scalar curvature $S_{g}$, there exists a positive
single peak solution $u_{\varepsilon}$ of (\ref{eq:main-sec2})
such that the peak point $q_{\varepsilon}$ approaches point $\bar{q}$
as $\varepsilon$ goes to zero.
\end{thm}

\begin{thm}
\label{thm:2-3}{\rm(See \cite{MP09b})} Assume that the scalar curvature
$S_{g}$ has $k\geq2$, $C^{1}$-stable critical points: $q_{1,}q_{2},\dots,q_{k}$.
Then, choosing $\varepsilon$ small enough, for any integer $j\le k$,
the problem (\ref{eq:main-sec2}) has a solution $u_{\varepsilon}$
with $j$ positive peaks $q_{1}^{\varepsilon},\dots,q_{j}^{\varepsilon}$
and $k-j$ negative peaks $q_{j+1}^{\varepsilon},\dots,q_{k}^{\varepsilon}$
such that $d_{g}(q_{i}^{\varepsilon},q_{i})\rightarrow0$ as $\varepsilon\rightarrow0$. 
\end{thm}
We now recall the definition of $C^{1}$-stable critical point.
\begin{defn}
Let $f\in C^{1}(M,\mathbb{R})$. The point $\bar{q}$ is a $C^{1}$-stable
critical point of $f$ if $\bar{q}$ is a critical point such that,
for any $\mu>0$ there exists $\delta$ for which any $h\in C^{1}(M,\mathbb{R})$
with\[
\max_{d_{g}(x,\bar{q})<\mu}\left[|f(x)-h(x)|+|\nabla_{g}f(x)-\nabla_{g}h(x)|\right]\leq\delta\]
has at least one critical point $q$ with $d_{g}(q,\bar{q})<\mu$.\end{defn}
\begin{thm}{\rm(See \cite{DMP09})} Let $\bar{q}\in M$ be an isolated local minimum
point of $S_{g}$. For each positive integer $k$, choosing $\varepsilon$
small enough, there exists a positive $k$-peaked solution $u_{\varepsilon}$
of (\ref{eq:main-sec2}) such that the $k$-peaks $q_{1}^{\varepsilon},\dots,q_{k}^{\varepsilon}$
collapse to $\bar{q}$, that is $d_{g}(q_{i}^{\varepsilon},\bar{q})\rightarrow0$
as $\varepsilon\rightarrow0$ for $i=1,\dots,k$. 
\end{thm}
By Theorem \ref{thm:2-1}, for a generic metric $g$, all the critical
points of $S_{g}$ are nondegenerate, then $C^{1}$-stable, isolated
and in a finite number. If $\nu$ is the number of critical points
of $S_{g}$, by Theorem \ref{thm:2-2} we have: $\nu$ positive solutions
with one peak, $\frac{\nu(\nu-1)}{2}$ solutions with two positive
peaks, $\dots$, one solution with $\nu$ positive peaks. Moreover,
by Theorem \ref{thm:2-3}, the problem (\ref{eq:main-sec2}) has some
sign changing solutions: for example $\frac{\nu(\nu-1)}{2}$ pairs
$(u_{\varepsilon},-u_{\varepsilon})$ of solutions with one positive
and one negative peak. Finally, since the global minimum point of
$S_{g}$ is isolated, the number of positive solutions of (\ref{eq:main-sec2})
goes to infinity as $\varepsilon$ goes to zero.

\section{Genericity of nondegeneracy for critical points for the mean curvature
of the boundary of a domain and applications\label{sec:3}}

In the following we denote by $E^{k}$ the vector space of all $C^{k}$
maps $\Psi:\mathbb{R}^{N}\rightarrow\mathbb{R}^{N}$ such that \[
\|\Psi\|_{k}=\sup_{x\in\mathbb{R}^{N}}\max_{\begin{array}{c}
0\le|\alpha|\le k\\
i=1,\dots,N\end{array}}\left|\frac{\partial^{\alpha}\Psi_{i}(x)}{\partial x_{1}^{\alpha_{1}}\cdots\partial x_{N}^{\alpha_{N}}}\right|<+\infty.\]

$E^{k}$ is a Banach space equipped with the norm $ $$ $$\|\cdot\|_{k}$.
Let $\mathscr{B}_{\rho}$ the ball in $E^{k}$ centered at zero with
radius $\rho$. Fixed an open bounded subset $\Omega\subset\mathbb{R}^{N}$
of class $C^{k}$, we have that the map \[
I+\Psi:\bar{\Omega}\rightarrow(I+\Psi)\bar{\Omega}\]
 is a diffeomorphism of class $C^{k}$ when $\Psi\in\mathscr{B}_{\rho}$
with $\rho$ small enough. We are interested in studying the nondegeneracy
of the critical points of the mean curvature of the boundary of the
domain $(I+\Psi)\Omega$ with respect to the parameter $\Psi$. More
precisely we can prove the following result (see \cite{MPpre})
\begin{thm}
\label{thm:3-1}Given a domain $\Omega\subset\mathbb{R}^{N}$ of class
$C^{k}$ with $k\ge3$, the set\[
A=\left\{ \begin{array}{c}
\Psi\in\mathscr{B}_{\rho}\subset E^{k}\ :\ \text{all the critical points of the mean curvature}\\
\text{of }(I+\Psi)\Omega\text{ are not degenerate}\end{array}\right\} \]
is an open dense subset of $\mathscr{B}_{\rho}$.
\end{thm}
We can apply Theorem \ref{thm:3-1} to study the following Neumann
problem\begin{equation}
\left\{ \begin{array}{ccc}
-\varepsilon^{2}\Delta_{g}u+u=|u|^{p-2}u &  & \text{in }\Omega\\
\frac{\partial u}{\partial\nu}=0 &  & \text{on }\partial\Omega\end{array}\right..\label{eq:neu-sec3}\end{equation}
Here $\Omega$ is a smooth bounded domain in $\mathbb{R}^{N}$, $N\geq2$
and $2<p<\frac{2N}{N-2}$ if $N\geq3$, $p>2$ if $N=2$. 

For Neumann problem (\ref{eq:neu-sec3}) we have some results about
existence of solutions in which the mean curvature $H$ of the boundary
$\partial\Omega$ of the domain $\Omega$ plays the same role of the
scalar curvature $S_{g}$ of the manifold $(M,g)$ in the problem
(\ref{eq:main-sec2}). We recall these results. Wei \cite{W1} and
Del Pino, Felmer, Wei \cite{DFW} proved that any $C^{1}$ stable
critical point of the mean curvature of $\partial\Omega$ gives a
single peaked solution. Gui \cite{G}, Li \cite{Li}, Wei and Winter
\cite{WW} proved the existence of multipeak solutions if the mean
curvature $H$ of $\partial\Omega$ has multiple $C^{1}$ stable critical
points. Dancer Yan \cite{DY} and Gui, Wei and Winter \cite{GWW}
proved the existence of clustered positive solutions such that the
peaks collapse to an isolated local minimum point of the mean curvature
of $\partial\Omega$. As far as it concerns the existence of sign
changing solutions there are results by Noussair and Wei \cite{NW98},
Micheletti and Pistoia \cite{MP08}, Wei and Weth \cite{WW05}, D'Aprile
and Pistoia \cite{DP11}. All these results require a sort of non
degeneracy of critical points of the mean curvature. By Theorem \ref{thm:3-1},
for a generic deformation $\Psi\in\mathscr{B}_{\rho}$, all the critical
points of the mean curvature of the boundary of the domain $(I+\Psi)\Omega$
are non degenerate, thus we can apply all the previous results.

\section{Generic properties for singularly perturbed nonlinear elliptic problems
on symmetric Riemannian manifolds \label{sec:4}}

We are now interested in studying the nondegeneracy of changing sign
solutions when the Riemannian manifold $(M,g)$ is symmetric. We consider
the problem \begin{equation}
\left\{ \begin{array}{ccc}
-\varepsilon^{2}\Delta_{g}u+u=|u|^{p-2}u &  & u\in H_{g}^{1}(M)\\
u(\tau x)=-u(x) &  & \forall x\in M\end{array}\right.\label{eq:ptau}\end{equation}
where $\tau:\mathbb{R}^{N}\rightarrow\mathbb{R}^{N}$ is an orthogonal
linear transformation such that $\tau\neq \Id$ and $\tau^{2}=\Id$,
$\Id$ being the identity on $\mathbb{R}^{N}$. Here the compact connected
Riemannian manifold $(M,g)$ of dimension $n\ge2$ is a regular submanifold
of $\mathbb{R}^{N}$ invariant with respect to $\tau$. In the following
we consider the space $H_{g}^{\tau}=\left\{ u\in H_{g}^{1}\ :\ \tau^{*}u=u\right\} $,
where the linear operator $\tau^{*}:H^{1}\rightarrow H^{1}$ is $\tau^{*}u=-u(\tau x)$. 

We obtain the following results (see \cite{GMp}) about the nondegeneracy
of changing sign solutions of (\ref{eq:ptau}) with respect the symmetric
metric $g$ considered as a parameter.
\begin{thm}
\label{thm:4-1}Given the metric $g_{0}$ on $M$ and the positive
number $\varepsilon_{0}$. The set\[
D=\left\{ \begin{array}{c}
h\in\mathscr{B}_{\rho}\ :\text{ any }u\in H_{g_{0}}^{\tau}\text{ solution of}\\
-\varepsilon_{0}^{2}\Delta_{g_{0}+h}u+u=|u|^{p-2}u\text{ is not degenerate}\end{array}\right\} \]
is a residual subset of $\mathscr{B}_{\rho}$. 
\end{thm}

\begin{thm}
\label{thm:4-2}Given the metric $g_{0}$ on $M$ and the positive
number $\varepsilon_{0}$. If there exists $\mu>m_{\varepsilon_{0},g_{0}}^{\tau}$
not a critical value of the functional $J_{\varepsilon_{0},g_{0}}$,
then the set\[
D=\left\{ \begin{array}{c}
h\in\mathscr{B}_{\rho}\ :\text{ any }u\in H_{g_{0}}^{\tau}\text{ solution of}\\
-\varepsilon_{0}^{2}\Delta_{g_{0}+h}u+u=|u|^{p-2}u\text{ is not degenerate}\end{array}\right\} \]
is an open subset of $\mathscr{B}_{\rho}$. 
\end{thm}
Here we set\begin{eqnarray*}
J_{\varepsilon_{0},g_{0}}(u) & = & \frac{1}{2}\int_{M}\left(\varepsilon_{0}^{2}|\nabla_{g_{0}}u|^{2}+u^{2}\right)d\mu_{g_{0}}-\frac{1}{p}\int_{M}|u|^{p}d\mu_{g_{0}}\\
\mathcal{N}_{\varepsilon_{0},g_{0}}^{\tau} & = & \left\{ u\in H_{g_{0}}^{\tau}\smallsetminus\left\{ 0\right\} \ :\ J'_{\varepsilon_{0},g_{0}}(u)\left[u\right]=0\right\} \\
m_{\varepsilon_{0},g_{0}}^{\tau} & = & \inf_{\mathcal{N}_{\varepsilon_{0},g_{0}}^{\tau}}J_{\varepsilon_{0},g_{0}}.\end{eqnarray*}
$\mathscr{B}_{\rho}$ is the ball centered at $0$ with radius $\rho$, with $\rho$ small
enough, in the Banach space $\mathscr{S}^{k}$, $k\ge$3 of all $C^{k}$ symmetric
covariant $2$-tensors $h(x)$ on $M$ such that $h(x)=h(\tau x)$
for $x\in M$ with $\|\cdot\|_{k}$ defined in (\ref{eq:normak}). 

If we choose the involution $\tau=-\Id$ and the symmetric manifold
such that $0\notin M$, using the previous theorems and the Morse theory
we obtain a generic result on the number of solutions which change sign exactly once.
\begin{thm}
Given the metric $g_{0}$ and $\varepsilon_{0}>0$ the set\begin{equation}
D=\left\{ \begin{array}{c}
h\in\mathscr{B}_{\rho}\ :\ \mbox{the equation }-\varepsilon_{0}^{2}\Delta_{g_{0}+h}u+u=|u|^{p-2}u\\
\mbox{has at least }P_{1}(M/G)\mbox{ pairs }(u,-u)\mbox{ of nontrivial solutions}\\
\mbox{which change sign exactly once}\end{array}\right\} \label{eq:D}\end{equation}
is a residual subset of $\mathscr{B}_{\rho}$
\end{thm}
Here $P_{t}(M/G)$ is the Poincar\`e polynomial of $M/G$ and $G=\{\Id, -\Id\}$
 
\begin{thm}
Given $g_{0}$ and $\varepsilon_{0}>0$, if there exists $\mu>m_{\varepsilon_{0},g_{0}}^{\tau}$
not a critical value of the functional $J_{\varepsilon_{0},g_{0}}$,
then the set $D$ defined in (\ref{eq:D}) is an open dense subset
of $\mathscr{B}_{\rho}$.
\end{thm}

\section{Appendix: A Transversality Theorem}

We recall transversality results which have been used to prove our genericity theorems. 
For the proofs we refer to \cite{Qu70,ST79,Uh76}.

\renewcommand{\labelenumi}{(\roman{enumi})}
\begin{thm}
\label{thm:gen1}Let $X,Y,Z$ be three real Banach spaces and let
$U\subset X,\ V\subset Y$ be two open subsets. Let $F$ be a $C^{1}$
map from $V\times U$ in to $Z$ such that 
\begin{enumerate}
\item For any $y\in V$, $F(y,\cdot):x\rightarrow F(y,x)$ is a Fredholm
map of index $0$.
\item $0$ is a regular value of $F$, that is $F'(y_{0},x_{0}):Y\times X\rightarrow Z$
is onto at any point $(y_{0},x_{0})$ such that $F(y_{0},x_{0})=0$.
\item The map $\pi\circ i:F^{-1}(0)\rightarrow Y$ is proper, where $i$
is the canonical embedding form $F^{-1}(0)$ into $Y\times X$ and
$\pi$ is the projection from $Y\times X$ onto Y
\end{enumerate}
Then the set \[
\theta=\left\{ y\in V\ :\ 0\text{ is a regular value of }F(y,\cdot)\right\} \]
 is a dense open subset of V
\end{thm}

\begin{thm}
\label{thm:gen2}If $F$ satisfies (i) and (ii) and
\begin{enumerate}
\addtocounter{enumi}{3}
\item The map $\pi\circ i$ is $\sigma$-proper,
that is $F^{-1}(0)=\cup_{s=1}^{+\infty}C_{s}$ where $C_{s}$ is a
closed set and the restriction $\pi\circ i_{|C_{s}}$ is proper for
any $s$ 
\end{enumerate}
then the set $\theta$ is a residual subset of $V$, that is $V\smallsetminus\theta$
is a countable union of closed subsets without interior points.
\end{thm}
\renewcommand{\labelenumi}{(\arabic{enumi})}

\end{document}